\documentclass{amsart}
\usepackage{eucal,amssymb}                                               
\begin{document}

\newcommand{\ad}{{\rm ad}}
\newcommand{\C}{{\mathbb C}}  
\newcommand{\Cee}{{\mathcal C}}
\newcommand{\cee}{{\mathfrak c}}
\newcommand{\cyc}{{\rm cyc}}                                
\newcommand{\D}{{\mathbb D}}                           
\newcommand{\e}{{\bf e}}
\newcommand{\ee}{{\exp_\infty}}
\newcommand{\ev}{{\rm ev}}
\newcommand{\EE}{{{\rm Exp}_\infty}}
\newcommand{\Ee}{{\bf E}}
\newcommand{\fr}{{\mathfrak {fr}}}  
\newcommand{\Faf}{{\mathfrak F}}
\newcommand{\faf}{{\mathfrak f}}
\newcommand{\G}{{\mathbb G}}                       
\newcommand{\Gal}{{\rm Gal}}
\newcommand{\GT}{{\rm GT}}
\newcommand{\grt}{{\mathfrak {grt}}}
\newcommand{\half}{{\textstyle{\frac{1}{2}}}}     
\newcommand{\K}{{\bf K}}                
\newcommand{\mot}{{\rm mot}}
\newcommand{\M}{{\mathcal M}}
\newcommand{\MM}{{\overline{\mathcal M}}}
\newcommand {\odd}{{\rm odd}}      
\newcommand{\Oh}{{\mathcal O}} 
\newcommand{\p}{{\mathfrak p}}                                       
\newcommand{\Q}{{\mathbb Q}}                      
\newcommand{\barQ}{{\overline{\mathbb Q}}}             
\newcommand{\MU}{{\bf MU}}
\newcommand{\R}{{\mathbb R}}
\newcommand{\Z}{{\mathbb Z}}

\newsymbol \Plh 207D

\title {The motivic Thom isomorphism} 
\author{Jack Morava}
\address{Department of Mathematics, Johns Hopkins University, Baltimore,
Maryland 21218}
\email{jack@math.jhu.edu}
\thanks{The author was supported in part by the NSF}
\date {15 November 2003}

\begin{abstract} 
The existence of a good theory of Thom isomorphisms in some rational category
of mixed Tate motives would permit a nice interpolation between ideas of 
Kontsevich on deformation quantization, and ideas of Connes and Kreimer on
a Galois theory of renormalization, mediated by Deligne's ideas on motivic
Galois groups. \end{abstract}

\maketitle
\section {Introduction} \bigskip

\noindent
This talk is in part a review of some recent developments in Galois theory, 
and in part conjectural; the latter component attempts to fit some ideas of 
Kontsevich about deformation quantization and motives into the framework of 
algebraic topology. I will argue the plausibility of the existence of liftings 
of  (the spectra representing) classical complex cobordism and $K$-theory 
to objects in some derived category of mixed motives over $\Q$. In itself 
this is probably a relatively minor technical question, but it seems to be 
remarkably consistent with the program of Connes, Kreimer, and others 
suggesting the existence of a Galois theory of renormalizations. \bigskip

\noindent     
{\bf 1.1} One place to start is the genus of complex-oriented manifolds 
associated to the Hirzebruch power series
\[
\frac{z}{\ee(z)}  =  z \Gamma(z)  =  \Gamma(1+z) 
\]
[25 \S 4.6]. Its corresponding one-dimensional formal group law is defined 
over the real numbers, with the entire function 
\[
\ee(z) = \Gamma(z)^{-1} : 0 \mapsto 0
\]
as its exponential. I propose to take seriously the related idea that the 
Gamma function
\[
\Gamma(z) \equiv z^{-1} \; {\rm mod} \; \R [[z]]
\]
defines some kind of universal asymptotic uniformizing parameter, or 
coordinate, at $\infty$ on the projective line, analogous to the role played
by the exponential at the unit for the multiplicative group, or the identity 
function at the unit for the additive group. \bigskip

\noindent
{\bf 1.2} The second  point of reference is a classical conjecture of Galois 
theory. The cyclotomic closure $\Q_\cyc$ of the rationals, defined by 
adjoining all roots of unity to $\Q$, is the maximal extension of $\Q$ with 
commutative Galois group; that group, isomorphic to the multiplicative
group ${\hat \Z}^\times$ of profinite integers, plays an important role in 
work of Quillen and Sullivan on the Adams conjecture and in differential 
topology. Shafarevich (cf. [27]) conjectures that the Galois group $\Gal
(\barQ/\Q_\cyc)$ is a {\bf free} profinite group; in other words, the full 
Galois group of $\barQ$ over $\Q$ fits in an exact sequence
\[
1 \to \Gal(\barQ/\Q_\cyc) \cong \widehat{\rm Free} \to \Gal(\barQ/\Q) \to 
\Gal(\Q_\cyc/\Q) \cong {\hat \Z}^\times \to 1 \;.
\]
What will be more relevant here is a related conjecture of Deligne [13 \S 
8.9.5], concerning a certain motivic analog of the Galois group which I will
denote $\Gal(\barQ/\Q)$, which is not a profinite but rather a proalgebraic 
groupscheme over $\Q$; it is in some sense a best approximation to the 
classical Galois group in this category, which should contain the original 
group as a Zariski-dense subobject. [I should say that calling this object 
a Galois group is an abuse of terminology; it is more properly described 
(cf. \S 4.4) as the motivic Tate Galois group of Spec $\Z$ (without reference 
to $\barQ$).] In any case, this motivic group fits in a 
similar extension
\[
1 \to \Faf_\odd  \to \Gal_\mot(\barQ/\Q) \to \G_m \to 1                 
\]
of groupschemes over $\Q$, where $\G_m$ is the multiplicative groupscheme, 
and  $\Faf_\odd$ is the prounipotent groupscheme defined by a free graded 
Lie algebra $\faf_\odd$ with one generator of each odd degree greater than 
one; the grading is specified by the action of the multiplicative group 
on the Lie algebra. The generators of this Lie algebra are thought to 
correspond with the odd zeta-values (via Hodge realization, [14 \S 2.14]) 
which are expected to be transcendental numbers (and thus outside the sphere 
of influence of Galois groups of the classical kind). \bigskip

\noindent
Kontsevich introduced his Gamma-genus in the context of the Duflo - Kirillov 
theorem in representation theory. He argued that it lies in the same orbit, 
under an action of the motivic Galois group, as the analog of the classical 
$\hat A$-genus [3 \S 8.5]. How this group fits in the topological context is 
less familiar, and to a certain extent this paper is nothing but an attempt 
to find a place for that group in algebraic topology. The history of this
question is intimately connected with Grothendieck's theory of anabelian 
geometry, and it enters Kontsevich's work through a conjectured Galois action 
on some form of the little disks operad. My impression these ideas are not 
yet very familiar to topologists, so I have included a very brief account of 
some of their history, with a few references, as an appendix below. \bigskip

\noindent
I should acknowledge here that Libgober and Hoffman [25,40] have studied a 
genus with related, but not identical, properties, and that an attempt to 
understand their work was instrumental in crystallizing the ideas behind this 
paper. I owe many mathematicians - including G. Carlsson, D. Christensen, 
F. Cohen, P. Deligne, A. Goncharov R. Jardine, T. Kohno, and T. Wenger - 
thanks for conversations about the material in this paper, and I am 
particularly indebted to a very knowledgeable and patient referee. In many 
cases they have saved me from mistakes and overstatements; but other such 
errors may remain, and those are the solely my responsibility. I also wish 
to thank the Newton Institute, and the Fields Institute program at Western 
Ontario, for support during the preparation of this paper. \bigskip

\section{The Gamma-genus} \bigskip
                                         
\noindent
{\bf 2.1} The Gamma-function is meromorphic, with simple poles at $z  = 0, -1,
-2, \dots$; we might therefore hope for a Weierstrass product of the form
\[
\Gamma(1+z)^{-1} \sim \prod_{n \geq 1} (1 + \frac{z}{n}) \;,
\]
from which we might hope to derive a power series expansion  
\[
\log \Gamma(1+z)  \sim - \sum_{n \geq 1} \log (1 + \frac{z}{n}) \sim 
\sum_{n,k \geq 1}(- \frac{z}{n})^k
\]
for its logarithm. Rearranging this carelessly leads to 
\[
\sum_{n,k \geq 1} \frac{(-z)^k}{k} \; \frac{1}{n^k} \sim \sum_{k \geq 1} 
\frac{\zeta(k)}{k} \; (-z)^k \;,
\] 
which is unfortunately implausible since $\zeta(1)$ diverges. In view of 
elementary renormalization theory, however, we should not be daunted: we 
can add `counter-terms' to conclude that 
\[
\log \prod_{n \geq 1} (1 +  \frac{z}{n}) \; e^{-z/n} \sim - \sum_{k \geq 2} 
\frac{\zeta(k)}{k} (-z)^k \;,
\]
and with a little more care we deduce the correct formula
\[
\Gamma(1+z) = \exp ( - \gamma z +  \sum_{k \geq 2} \frac{\zeta(k)}{k} \; 
(-z)^k) \;,
\]
where $\gamma$ is Euler's constant. Reservations about the logic of this 
argument may perhaps be dispelled by observing that  
\[
\Gamma(1+z) \; \Gamma(1-z)  =  \exp ( \sum_{k \geq 1} \frac{\zeta(2k)}{k} \; 
z^{2k}) \;;
\]
Euler's duplication formula implies that the left-hand side equals 
\[                                                     
z \Gamma(z) \; \Gamma(1-z) = \frac{\pi z}{\sin \pi z} \;,
\]
consistent with the familiar evaluation of $\zeta$ at positive even integers 
in terms of Bernoulli numbers. \bigskip

\noindent
{\bf 2.2} From this perspective, the Hirzebruch series 
\[
\Gamma(1+z) = (\frac{\pi z}{\sin \pi z})^\half \; \exp ( - \gamma z + \sum 
\frac {\zeta(\odd)}{\odd} \; z^\odd ) \;;
\]
for Kontsevich's genus does in fact look like some kind of deformation of 
the $\hat A$-genus; its values on a complex-oriented manifold will be 
polynomials in odd zeta-values, with rational coefficients. Similarly, the 
Witten genus 
\[
\phi_W(x) = \frac{x/2}{\sinh x/2} \; \prod_{n \geq 1} [(1 - q^n u)
(1 - q^nu^{-1})]^{-1}
\]
([57], with $u = e^x$) can be written in the form
\[                                      
\exp ( - 2 \sum_{k \geq 1} g_k \frac{x^k}{k!}) \;,
\]
where the coefficients $g_k$ are modular forms, with $g_\odd = 0$: it is 
also a deformation of $\hat A$, in another direction. \bigskip

\noindent
[Behind the apparent discrepancies in these formulae is the issue of complex 
versus oriented cobordism: there are several possible conventions relating 
Chern and Pontrjagin classes. Hirzebruch expresses the latter as symmetric 
functions of indeterminates $x_i^2$, and writes the genus associated to the 
formal series $Q(z)$ as $\prod Q(x_i^2)$; thus for the $\hat A$-genus, 
\[
Q(z) = \frac {\half \surd z}{\sinh \half \surd z} \;.
\]
An alternate convention, used here, writes this symmetric function in the form
\[
\prod ((\frac{x_i/2}{\sinh x_i/2})^\half \cdot (\frac{-x_i/2}{\sinh (- x_i/2)}
)^\half)  \;.
\]
The relation between the indeterminates $x$ and $z$ is a separate issue; I 
take $z$ to be $2 \pi i x$.] \bigskip

\noindent
Kontsevich suggests that the values of the zeta function at odd positive 
integers (expected to be transcendental) are subject to an action of the 
motivic group $\Gal_\mot(\barQ/\Q)$, and that the $\hat A$-genus 
and his $\Gamma$ genus lie in the same orbit of this action. One natural way 
to understand this is to seek an action of that group on genera, and 
thus on the complex cobordism ring; or, perhaps more naturally, on some form 
of its representing spectrum. Before confronting this question, it may be 
useful to present a little more background on these zeta-values. \bigskip

\section{Symmetric and quasisymmetric functions}
\bigskip                                                    

\noindent
{\bf 3.1} The formula
\[
\prod_{k \geq 1}(1 + x_kz) = \sum_{k \geq 0} e_k z^k  =  \exp( - \sum_{k 
\geq 1} \frac{p_k}{k} (-z)^k) \;,
\]
where $e_k$ is the $k$th elementary symmetric function, and $p_k$ is the 
$k$th power sum, can be derived by formal manipulations very much like those 
in the preceding section, by expanding the logarithm of $\prod (1 + x_n z)$;
such arguments go back to Newton. The specialization 
\[
x_k \mapsto k^{-2}
\]                            
(cf. the second edition of MacDonald's book [43 Ch I \S 2 ex 21]) leads to 
Bernoulli numbers, but the map                    
\[
x_k \mapsto k^{-1}
\]
is trickier, because of convergence problems like those mentioned above; it 
defines a homomorphism from the ring of symmetric functions to the reals, 
sending $p_k$ to $\zeta(k)$ when  $k > 1$, while $p_1 \mapsto \gamma$ [24]. 
Under this homomorphism the even power sums $p_{2k}$ take values in the 
field $\Q(\pi)$. \bigskip
                                                       
\noindent
{\bf 3.2} The Gamma-genus is thus a specialization of the formal group law 
with exponential
\[
\EE(z) = \frac{z}{e(z)} = z \prod_{k \geq 1}(1 + x_k z)^{-1} = \sum_{k \geq 0}
(-1)^k h_k z^{k+1}
\]
having the complete symmetric functions (up to signs) as its coefficients. 
This is a group law of {\bf additive} type: its exponential, and hence its 
logarithm, are both defined over the ring of polynomials generated by 
the elements $h_k$. This group law is classical: it is defined by the 
Boardman-Hurewicz-Quillen complete Chern number homomorphism
\[
MU^*(X) \to H^*(X,\Z[h_*])
\]
defined on coefficients by the homomorphism
\[
{\bf Lazard} \to {\bf Symm}
\]
from Lazard's ring which classifies the universal group law of additive type. 
The Landweber-Novikov Hopf algebra $S_* = \Z[t_*]$ represents the prounipotent
groupscheme $\D_0$ of formal diffeomorphisms 
\[
z \mapsto t(z) = z + \sum_{k \geq 1} t_k z^{k+1}
\]
of the line, with coproduct
\[
\Delta (t(z)) = (t \otimes 1)((1 \otimes t)(z)) \in (S_* \otimes S_*)[[z]] \;.
\]
The universal group law $t^{-1}(t(X) + t(Y))$ of additive type is thus 
classified by the homomorphism
\[
{\bf Lazard} \to {\bf Lazard} \otimes S_* \to \Z \otimes S_* \to {\bf Symm}
\]
representing the orbit of the Thom map ${\bf Lazard} \to \Z$ (which classifies 
the additive group law) under the action of $\D_0$. This identifies the 
algebra $S_*$ with the ring of symmetric functions by $t_k \mapsto (-1)^kh_k$. 
\bigskip

\noindent
{\bf 3.3} The symmetric functions are a subring 
\[
{\bf Symm} \to {\bf QSymm}
\]
of the larger ring of {\bf quasi}symmetric functions, which is Hopf dual to 
the universal enveloping algebra $\Z \langle Z_1,\dots \rangle$ of the free 
graded Lie algebra $\faf_*$ with one generator in each positive degree [23], 
given the cocommutative coproduct
\[
\Delta Z_i = \sum_{i = j+k} Z_j \otimes Z_k \;.
\]
Standard monomial basis elements for this dual Hopf algebra, under 
specializations like those discussed above [6 \S 2.4, 23], map to polyzeta 
values
\[
\zeta (i_1,\dots,i_k) = \sum_{n_1 > \dots > n_k \geq 1} \frac{1}{n_1^{i_1} 
\cdots n_k^{i_k}} \in \R \;;
\] 
note that there are convergence difficulties unless $i_1 > 1$. If we think of 
the Gamma-genus as taking values in the field $\Q(\zeta) \subset \R$ generated
by such polyzeta values, then it is the specialization of a homomorphism
\[
{\bf Lazard} \to {\bf Symm} \to {\bf QSymm}
\]
representing a morphism from the prounipotent groupscheme $\Faf$ with Lie 
algebra $\faf_*$ to the moduli space of one-dimensional formal group laws. 
\bigskip

\noindent
Because we are dealing with group laws of additive type, there seems to be 
little loss in working systematically over a field of characteristic zero, 
where Lie-theoretic methods are available. Over such a field any formal group 
is of additive type: the localization of the map from the Lazard ring to the 
symmetric functions is an isomorphism. Similarly, over the rationals the 
Landweber - Novikov algebra is dual to the enveloping algebra of the Lie 
algebra of vector fields
\[
z_k = z^{k+1} \partial/\partial z \;, k \geq 1
\]
on the line, and the embedding of the symmetric in the quasisymmetric 
functions sends the free generators $Z_k$ to the Virasoro generators $z_k$, 
corresponding to a group homomorphism $\Faf \to \D_0$. \bigskip

\noindent
{\bf 3.4} It will be useful to summarize a few facts about Malcev completions
and pro-unipotent groups [13 \S 9]. The rational group ring $\Q[G]$ of a 
discrete group $G$ has a natural lower central series filtration; its 
completion $\widehat{Q[\G]}$ with respect to that filtration is a topological 
Hopf algebra, whose continuous dual represents a pro-unipotent groupscheme 
over $\Q$. Applied to a finitely-generated free group, for example, this 
construction yields a Magnus algebra of noncommutative formal power series. 
\bigskip

\noindent
There are many variations on this theme: in particular, an action of the 
multiplicative group $\G_m$ defines a grading on a Lie algebra. The action 
\[
t_k \mapsto u^k t_k \;,                                          
\]
($u$ a unit) on the group of formal diffeomorphisms defines an extension of 
its Lie algebra by a new Virasoro generator $v_0$, corresponding to an 
extension $S_*[t_0^{\pm}]$ of the Landweber-Novikov algebra. The group of 
formal diffeomorphisms is pro-unipotent, and this enlarged object is most 
naturally interpreted as a semidirect product $\D_0 \rtimes \G_m$. Grading 
the free Lie algebra $\faf$ similarly extends the homomorphism above to 
\[
\Faf \rtimes \G_m \to \D \rtimes \G_m \;.
\]

\section{Motivic versions of classical $K$-theory and cobordism}
\bigskip

\noindent
{\bf 4.1} There are now several (eg [39, 55]) good and probably equivalent 
constructions of a triangulated category $DM(k)$ of motives over a field $k$ 
of characteristic zero. The subject is deep and fascinating, and I know at 
best some of its vague outlines. Since this paper is mostly inspirational,
I will not try to provide an account of that category; but as it is after all 
modelled on spectra, it is perhaps not too much of a reach to think that 
some of its aspects will look familiar to topologists. \bigskip

\noindent
One approach to defining a motivic category $DM(k)$ starts from a category 
whose morphisms are elements of a group of algebraic correspondences. At 
some later point it becomes useful to tensor these groups with $\Q$, 
resulting in a category $DM_\Q(k)$ whose Hom-objects are rational vector 
spaces. The underlying concern of this paper is the relation of such motivic 
categories to classical topology; but stable homotopy theory over the 
rationals is equivalent to the theory of graded vector spaces. This has 
the advantage of rendering some of the conjectures below almost trivially 
true -- and the disadvantage of making them essentially contentless. 
Behind these conjectures, however, lies the hope that they might say
something {\bf before} rationalization, and for that reason I have outlined 
here a rough theory of {\bf integral} geometric realizations of motives: 
\bigskip

\noindent
{\bf 4.2} The category $DM(k)$ contains certain canonical Tate objects 
$\Z(n)$, defined [55 \S 2.1, but see also 13 \S 2] as tensor powers of a 
reduced version $\Z(1)$ of the projective line. Grothendieck's original 
category of `pure' motives, constructed from smooth projective varieties, 
is (in some generality [28]) semisimple, but categories of motives built 
from more general (non-closed) varieties admit nontrivial extensions. The 
(derived) category $DMT_\Q(k)$ of mixed Tate motives can be defined as the 
smallest tensor triangulated subcategory of $DM_\Q(k)$ containing the Tate 
objects. In this rationalized category, it is natural to denote 
the (images of) the generating objects by $\Q(n)$; however, I will be most
interested here in the case $k = \Q$ and in a certain more subtle construction
of a (rationalized, though I will now drop the subscript) category $DMT(\Z)$
of $DMT_\Q(\Q)$ [14 \S 1.6], closely related to the motives over Spec $\Z$ 
`with integral coefficients' in the sense of [13 \S 1.23, 2.1]. This is still 
a $\Q$-linear category, but its objects have stronger integrality properties 
than one might naively expect. \bigskip

\noindent
In particular: one of the foundation-stones of the theory of mixed motives 
is an isomorphism
\[
{\rm Ext}^1_{MT(\Z)}(\Q(0),\Q(n)) \cong K_{2n-1}(\Z) \otimes \Q
\]
(cf. [2, 13 \S 8.2]. As the referee points out, one has to be careful
here; the corresponding description for the category $MT(\Q)$ involves the
algebraic $K$-theory of $\Q$, which is much larger than that of $\Z$).
The groups on the right have rank one for {\bf odd} $n > 1$, and vanish 
otherwise, by work of Borel; the theory of regulators says that to some 
extent the zeta-values 
\[
\frac{(n-1)!}{(2 \pi i)^n} \; \zeta(n)
\]
(cf. [13 \S 3.7]) can be interpreted as natural generators for these 
groups. This is strikingly reminiscent to a homotopy-theorist of the 
identification (for $n$ {\bf even}) of the
group
\[
{\rm Ext^1}_{\rm Adams}(K(S^0),K(S^{2n})) 
\]
of extensions of modules over the Adams operations, with the cyclic subgroup
of $\Q/\Z$ generated by this zeta-number. These connections 
between the image of the $J$-homomorphism and the groups $K_{4k-1}(\Z)$, 
go back to the earliest days of algebraic $K$-theory [16, 50]. \bigskip

\noindent
{\bf 4.3} This suggests that there might be some use for a notion of 
geometric or homotopy-theoretic realization for motives, which manages to 
retain some integral information. Aside from tradition (algebraic geometers 
usually work with cycles over $\Q$, and topologists have been neglecting 
correspondences since Lefschetz), there seems to be no obstacle to the 
development of such a theory. Indeed, let $E$ be a multiplicative (co)homology
functor (ie a ring-spectrum), supplied with a natural class of $E$-orientable 
manifolds: if $E$ were stable homotopy, for example, we could use stably 
parallelizeable manifolds. In the case of interest below, however, $E$ will 
be complex $K$-theory, and the manifolds will be smooth (proper) algebraic 
varieties over $\C$. \bigskip

\noindent
Such manifolds ($X,Y,\dots$) define an additive category ${\rm Corr}_E$
with $E_*(X \times Y)$ (suitably graded) as Hom-objects; composition of such
morphisms can be defined using Pontrjagin-Thom transfers [45]. It is 
straightforward to check that 
\[
X \mapsto X_+ \wedge \Ee : {\rm Corr}_E \to ({\rm Spectra})
\]
is a functor: the necessary homomorphism
\[
E_*(X \times Y) = [S^*, X_+ \wedge Y_+ \wedge \Ee] \to [X_+ \wedge \Ee,Y_+ 
\wedge \Ee]_* 
\]
of Hom-objects is defined by the adjoint composition
\[
X_+ \wedge Y_+ \wedge \Ee \wedge X_+ \wedge \Ee \cong X_+ \wedge X_+ \wedge 
Y_+ \wedge \Ee \wedge \Ee \to Y_+ \wedge \Ee
\]
built from the multiplication map of $\Ee$ and the composition
\[
X_+ \wedge X_+ \wedge \Ee \to X_+ \wedge \Ee \to \Ee
\]
of the transfer $\Delta_!$ associated to the diagonal map, with the projection
of $X$ to a point. \bigskip

\noindent
Following the pattern laid out by Voevodsky, we can now define a category 
of (topological) `$E$-motives', and when $E = K$ it is a classical fact 
[1] that an algebraic cycle defines a nice $K$-theory class. This allows 
us to associate to an embedding of $k$ in $\C$, a triangulated `realization' 
functor 
\[
\K_\mot : X \mapsto X(\C)_+ \wedge \K : DM(k) \to ({\rm Spectra}) \;.
\]
\medskip

\noindent
{\bf 4.4} When $k$ is a number field, $DMT_\Q(k)$ possesses a theory of 
truncations, or $t$-structures [14, 37], analogous to the Postnikov systems 
of homotopy theory; the {\bf heart} of this structure is an abelian 
tensor category $MT_\Q(k)$ of mixed Tate motives. Its existence permits us 
to think of $DMT_\Q(k)$ as the derived category of $MT_\Q(k)$; in particular, 
the (co)homology of an object of the larger category becomes in a natural way 
[31 \S 2.4] an object of $MT_\Q(k)$. Similar considerations hold for the
more rigid category $DMT(\Z)$, and since we are working in a rational, stable 
context, I will write $\pi_*$ for the homology groups of an object in
this category, given this enriched structure. [For the purposes of this
presentation I've reversed the logical order of construction: in fact in
[14] the category $MT(\Z)$ is constructed first.] \bigskip

\noindent
Now under very general conditions (involving a suitably rigid duality), an 
abelian tensor category with rational Hom-objects can be identified with a 
category of representations of a certain groupscheme of automorphisms of a 
suitable forgetful functor on the category; the resulting groupscheme is 
called a motivic (Galois) group. This theory applies to $MT(\Z)$, and as
was noted in the introduction, $\Gal_\mot (\barQ/\Q)$ is the corresponding 
groupscheme [13 \S 8.9.5, 14 \S 2; cf. also 2 \S 5.10]. In the preceding 
paragraph we constructed a homological functor 
\[
K_\mot := \pi_* \K_\mot
\]
from $DM_\Q$ to $\Q$-vector spaces, together with a preferred lift of the
functor to the category of spectra. It is easy to see that
\[
\pi_\odd  \K_\mot = 0 \; , \; {\rm while} \; \pi_{2n}\K_\mot = \Q(n) \;,
\]
as representations of $\Gal_\mot(\barQ/\Q)$, and thus that $K_\mot$ is 
represented by the mixed (Bott!)-Tate object 
\[
\oplus_{n \in \Z}\Q(n)[2n] = \Q[b^\pm] \;;
\]
in other words, we have constructed a lifting of the rationalized classical
$K$-theory functor to the category of mixed Tate motives, with an action of
$\Gal_\mot(\barQ/\Q)$ which factors through the multiplicative quotient.
\bigskip

\noindent
The possible existence of a descent spectral sequence for the automorphisms
of the $K$-theoretic realization functor of \S 4.3 seems to be an interesting
question, especially when restricted to some category of mixed Tate motives.
\bigskip

\noindent
{\bf 4.5} The {\bf main conjecture} of this paper is that, similarly, a 
rational version of complex cobordism lifts to an object $\MU_\mot 
\in DMT(\Z)$, with an action of $\Gal_\mot(\barQ/\Q)$ on $\pi_* \MU_\mot$ 
defined by the obvious embedding 
\[
\Faf_\odd \rtimes \G_m \to \Faf \rtimes \G_m
\]
followed by a homomorphism from the latter group to the diffeomorphisms of 
the formal line, cf. \S 2.4 above. \bigskip

\noindent
This is a conjecture about an object characterized by its universal 
properties, so it can be reformulated in terms of the structures thus 
classified. The theory of Chern classes is founded on Grothendieck's 
calculation of the cohomology of the projectification $P(V)$ of a vector 
bundle $V$ over a scheme $X$. It follows immediately from his result that 
the cohomology of Atiyah's model 
\[
X^V := P(V \oplus 1) /P(V) 
\]
for a Thom space as a relative motive is free on one generator over that of 
$X$. Such a generator is a Thom class for $V$, but in the motivic context 
there seems to be no natural way to construct such a thing; this is related
to the inconvenient nonexistence of abundantly many sections of vector 
bundles in the algebraic category. \bigskip

\noindent
For a systematic theory of Thom classes it is enough, according to the 
splitting principle, to work with line bundles $L$, and in this context it is 
relevant that the Thom complexes
\[
X^L = P(L \oplus 1)/P(L) \cong P(1 \oplus L^{-1})/P(L^{-1}) = X^{L^{-1}}
\]
of a line bundle and its reciprocal are isomorphic objects. The conjecture 
about $\MU_\mot$ can be thus reformulated in terms of a theory of motivic 
Thom and Euler classes $U_\mot(L), \e(L)$ for line bundles $L$, satisfying 
a {\bf motivic Thom axiom}
\[
U_\mot(L^{-1}) = - U_\mot(L) \;,\; \e(L^{-1}) = - \e(L) \;;
\]
the conjecture is then the assertion that an element $\sigma \in 
\Gal_\mot(\barQ/\Q)$ sends $U_\mot(L)$ to another Thom class 
\[
\sigma (U_\mot(L)) = [1 + \sum_{k > 0}  \sigma_k \e(L)^k] \cdot U_\mot(L)
\]
for $L$, with coefficients $\sigma_k$ depending only on $\sigma$. Since  
\[                                                     
\sigma(U_\mot(L^{-1})) =  - \sigma(U_\mot(L)) \;,
\]
it follows that
\[
[1 + \sum_{k > 0} \sigma_k (- \e(L))^k] \cdot (- U_\mot(L)) = - [1 + \sum_{k >
0} \sigma_k \e(L)^k] \cdot U_\mot(L) \;,
\]
which entails that the classes $\sigma_\odd = 0$, distinguishing the Hopf 
subalgebra   
\[
S_\ev = \Z[t_{2k} \; | \; k > 0]
\]
which represents the group of {\bf odd} diffeomorphisms of the formal line 
[46 \S 3.3]. Away from the prime two, classical complex cobordism is a 
kind of base extension
\[
{\bf MU}[1/2] \sim {\rm SO/SU} \wedge {\bf MSO}[1/2]
\]
of oriented cobordism, and I'm suggesting the existence of a similar
splitting for the hypothetical motivic lift of complex cobordism. \bigskip 

\noindent
{\bf 4.6} After this paper had been submitted for publication, I became 
aware of the very elegant recent work of Levine and Morel [38], where
an algebraic cobordism functor is characterized as a universal cohomology
theory on the category of schemes, endowed with pullback and pushforward 
transformations satisfying certain natural axioms pf compatibility. 
[Voevodsky [56] has also considered a motivic version of the cobordism 
spectrum; its relation with their work is discussed briefly in the 
introduction to their paper.] I believe their work is fundamentally
compatible with the conjectures made here, given a slight difference in
framework and emphasis: they suppose a Thom isomorphism (or, equivalently,
a system of covariant transfers) is to be given as part of the
structure of a cohomology theory on schemes, while the spectrum hypothesized 
here is merely a ringspectrum, with no preferred choice of orientation.
\bigskip

\noindent
{\bf 4.7} I should note that the trivial action of $\Faf_\odd$ on 
$\pi_*(\K_\mot)$, together with the usual action of $\G_m$ defined by the 
grading, is consistent with the existence of a spectral sequence
\[
E_2^{*,*} = H_c^*(\Gal_\mot(\barQ/\Q), K_\mot^*) \Longrightarrow K^*(\Z) 
\otimes \Q 
\]
of descent type: using the Hochschild-Serre spectral sequence for the 
semidirect product decomposition of the Galois group (and confusing continuous
cochain with Lie algebra cohomology) we start with
\[
H^*_c(\Faf_\odd,K_\mot^*) \cong \Q \langle e_{2k+1} \; | \; k \geq 1 \rangle 
[b] \;,
\]
where angled brackets denote a vector space spanned by the indicated elements 
(with $e_{2r+1}$ in degree $(1,0)$, and the Bott (-Tate?!) element $b$ in 
degree $(0,-2)$). The $\G_m$-action sends $b$ to $ub$, where $u$ is a unit 
in whatever ring we're over; similarly, 
\[
e_{2k+1} \mapsto u^{-2k-1} e_{2k+1} \;.
\]
Thus $e_{2k+1} b^{2k+1} \in E_2^{1,-4k-2}$ is $\G_m$-invariant, yielding a 
candidate for the standard generator in  $K_{4k+1}(\Z) \otimes \Q$. \bigskip
                              
\section{Quantization and asymptotic expansions}
\bigskip
                                                       
\noindent
{\bf 5.1} The motivic Galois group appears in Kontsevich's work through a 
conjectured action on deformation quantizations of Poisson structures [cf.
[53]]. The framework of this paper suggests a plausibly related action 
on an algebra of asymptotic expansions for geometrically defined functionals 
on manifolds, interpreted in terms of the cobordism ring of symplectic 
manifolds [18, 19]. This is isomorphic to the complex cobordism (abelian Hopf)
algebra $MU_*(B\G_m(\C))$ of circle bundles, and the dual (rationalized) Hopf 
algebra 
\[
MU^*_\Q(B\G_m(\C)) \cong MU^*_\Q[[\Plh]] \;,
\]                                      
where $\Plh = \sum_{k \geq 1} \C P_{k-1} \e^k/k$, can be interpreted [46] 
as an algebra generated by the coefficients of a kind of universal asymptotic 
expansion for geometrically defined heat kernels (or Feynman measures, via 
the Feynman-Kac formula [20 \S 3.2]), as the Chern class $\e$ of the circle 
bundle approaches infinity. \bigskip

\noindent
A Poisson structure on an even-dimensional manifold $V$ is a bivector field 
(a section of the bundle $\Lambda^2 T_V$) satisfying a Jacobi identity 
modelled on that satisfied by the inverse of a symplectic structure. A 
symplectic manifold is thus Poisson, and although I am aware of no useful 
notion of Poisson cobordism (but cf. [7]) one expects a natural restriction 
map from asymptotic invariants Poisson manifolds to the corresponding ring 
for symplectic manifolds. If the conjectures above are correct, then one 
might further hope that (some motivic version of) such a restriction map 
would be equivariant, with respect to some motivic Galois action. \bigskip

\noindent
{\bf 5.2} Working in the opposite (local to global) direction, Connes and 
Kreimer have recently developed a systematic program for understanding 
classical quantum field theoretic renormalization in terms of its symmetries. 
The standard methods (eg dimensional regularization) for dealing with the
singular integrals which appear in classical perturbation theory replaces 
them with certain meromorphic functions, and through work of Broadhurst,
Kreimer, and others it has become more and more clear that the polar 
coefficients of these meromorphic functions are frequently elements of 
the polyzeta algebra. [Kontsevich has suggested that this is always so, but 
his program of proof fails, by the arguments of [4], and at present the 
question seems to be open.] Connes and Kreimer[12, 21, 26, 35 \S 12] have
developed a systematic approach to the theory of Feynman integrals through 
certain Hopf algebras related to automorphism groups [42] of operads defined 
by graphs of various sorts. \bigskip

\noindent
There are deep connections between the Grothendieck - Teichm\"uller group 
and the Lie algebras of these automorphism groups [29, 54], and it 
seems likely that they (and the theory of quasisymmetric functions, via free 
Lie algebras) will eventually be understood to be intimately related; the 
appearance of polyzeta values in the theory of quantum knot invariants (cf. 
eg [36]) is another source of recent interest in this subject. Perhaps the
deepest (and most precise) approach to the relations between these topics
may be the work of Goncharov, who associates to a field $F$ a certain
Hopf algebra ${\mathcal T}_\bullet(F)$ of $F$-decorated planar trivalent
trees and a closely related Hopf algebra ${\mathcal I}_\bullet(F)$ of motivic
iterated integrals. According to the correspondence principle of [21 \S 7],
the renormalization Hopf algebra corresponding to certain types of Feynman
integrals should be closely related to a precisely defined subgroup of the
motivic Tate Galois group of $F$. A slight strengthening of this 
correspondence principle would settle the question of the role of polyzeta 
values in perturbative expansions of Feynman integrals. \bigskip

\noindent
{\bf 5.3} In exemplary cases Connes and Kreimer construct a very interesting 
representation of the prounipotent groupscheme underlying their 
renormalization algebra in the group of {\bf odd} formal diffeomorphisms of 
the line [10 \S 1 eq. 20, \S 4 eq. 2]. It also seems quite possible that the 
action of their groupscheme on asymptotic expansions defined by Feynman 
measures associated to suitable Lagrangians [30] factor through an 
action of the motivic Galois group on cobordism, along the lines suggested 
in \S 4.5 above. \bigskip                       
                                                  
\noindent
{\bf appendix: motivic models for the little disk operad}
\bigskip \bigskip

\noindent
{\bf 1} In 1984 Grothendieck suggested the study of the action of the Galois 
group $\Gal(\barQ/\Q)$ as automorphisms of the moduli of algebraic curves, 
understood as a collection of stacks linked by morphisms representing various 
geometrically natural fusion operations. There are remarkable analogies 
between his ideas and contemporary work in physics on conformal field 
theories, and in 1990 Drinfel'd [15] unified at least some of these lines of 
thought by constructing a pronilpotent group $\GT$ of automorphisms of 
certain braided tensor categories, together with a faithful representation 
of the absolute Galois group in that group. \bigskip

\noindent
This program has been enormously productive; the LMS notes [42, 51] are one 
possible introduction to this area of research, but evolution has been 
extremely rapid. In the late 90's Kontsevich [33, 34] recognized connections 
between these ideas, Deligne's question on Hochschild homology, deformation 
quantization, and other topics, while physicists [9-12] interested in the 
algebra of renormalization were developing sophisticated Hopf-algebraic
techniques, which are now believed [5 \S 8,9] to be closely related to the 
Hopf-algebraic constructions of Drinfel'd. \bigskip

\noindent
The point of this appendix is to draw attention to a central conjecture in 
this circle of ideas: that the Lie algebra of $\GT$, which acts as 
automorphisms of the system of Malcev completions of the braid groups, is a 
{\bf free} graded Lie algebra, with one generator in each odd degree greater
than one. The braid groups in fact form an operad, and I want to propose the 
related problem of identifying the automorphisms of the operad of Lie 
algebras defined by the braid groups (cf. [8]), in hope that this will shed 
some light on this question, and the closely related conjecture that 
Deligne's motivic group acts {\bf faithfully} on the unipotent motivic 
fundamental group [13] of ${\mathbb P)}_1 -\{0,1,\infty\}$ (with nice 
tangential base point). \bigskip

\noindent
{\bf 2} For the record, an operad (in some reasonable category) is a 
collection of objects $\{\Oh_n,\; n \geq 2 \}$ together with composition 
morphisms
\[
c_I : \Oh_{r(I)} \times \prod_{i \in I} \Oh_{i} \to \Oh_{|I|}
\]
where $I = i_1,\dots,i_r$ is an ordered partition of $|I| = \sum i_k$ with 
$r(I)$ parts. These compositions are subject to a generalized associativity 
axiom, which I won't try to write out here; moreover, the operads in this 
note will be {\bf permutative}, which entails the existence of an action of 
the symmetric group $\Sigma_n$ on $\Oh_n$, also subject to unspecified axioms.
Not all of the operads below will be unital, so I haven't assumed the 
existence of an object $\Oh_1$; but in that case, and under some mild 
assumptions, an operad can be described as a monoid in a category of objects 
with symmetric group action, with respect to a somewhat unintuitive product 
[cf. eg. [17]]. \bigskip

\noindent
The moduli $\{ \MM_{0,n+1} \}$ of stable genus zero algebraic curves marked
with $n+1$ ordered smooth points form such an operad: if the final marked 
point is placed at infinity, then composition morphisms are defined by gluing 
the points at infinity of a set of $r$ marked curves to the marked points 
{\bf away} from infinity on some curve marked with $r+1$ points. This is an
operad in the category of algebraic stacks defined over $\Z$, so the Galois 
group $\Gal(\barQ/\Q)$ acts on many of its topological and cohomological 
invariants. \bigskip

\noindent
By definition, a stable algebraic curve possesses at worst ordinary double 
point singularities. The moduli of {\bf smooth} genus zero curves thus define 
a system $\M_{0,*} \subset \MM_{0,*}$ of subvarieties, but not a suboperad: 
the composition maps glue curves together at smooth points, creating new 
curves with nodes from curves without them. However, the spaces $\M_{0,n+1}
(\C)$ of complex points of these moduli objects have the homotopy types of the
spaces of configurations of $n$ distinct points on the complex line, which are
homotopy-equivalent to the spaces of the little disks operad  $\Cee_2$. 
\bigskip

\noindent
Behind the conjectures of Deligne, Drinfel'd, and Kontsevich lies an 
apparently unarticulated question: {\bf is the little disks operad defined 
over $\Q$?}; or, more precisely: is there a version of the little disks 
operad in which the morphisms, as well as the objects, lie in the category 
of algebraic varieties defined over the rationals? Thus in [33] (end of \S 
4.4) we have \bigskip

\begin{quote}

The group $\GT$ maps to the group of automorphisms in the homotopy sense of 
the operad Chains($\Cee_2$). Moreover, it seems to coincide with 
Aut(Chains($\Cee_2$)) when this operad is considered as an operad not of 
complexes but of differential graded cocommutative coassociative coalgebras 
\dots

\end{quote} \bigskip

\noindent
and in [34] (end of \S 3) \bigskip

\begin{quote}

There is a natural action of the Grothendieck-Teichm\"uller groups on the 
rational homotopy type of the [Fulton-MacPherson version of the] little disks 
operad \dots

\end{quote} \bigskip

\noindent
although the construction in [34 \S 7] is not (apparently) defined by 
algebraic varieties. It may be that the question above is naive 
[cf. [53]], but a positive answer would imply the existence of a system 
of homotopy types with action of the Galois groups, whose algebras of chains 
would have the properties claimed above; moreover, the system of fundamental 
groups of these homotopy types would yield an action of the Galois group on 
the system of braid groups, suitably completed. \bigskip

\noindent
{\bf 3} This is probably a good place to note that $\Gal(\barQ/\Q)$ acts by 
automorphisms of the etal\'e homotopy type of a variety defined over $\Q$, 
which is not at all the same as a continuous action on the space of complex 
points of the variety. In fact one expects to recover classical invariants 
of a variety (the cohomology, or fundamental group, of its complex points, 
for example) only up to some kind of completion. Various kinds of invariants 
[etal\'e, motivic, Hodge-deRham, \dots] each have their own associated 
completions, some of which are still quite mysterious; the fundamental group, 
in particular, comes in profinite [49] and prounipotent versions. \bigskip

\noindent
Drinfel'd works with the latter, which corresponds to the Malcev completion 
used in rational homotopy theory. A free group on $n$ generators corresponds 
to the graded Lie algebra defined by $n$ noncommuting polynomial generators, 
and the Lie algebra $\p_n$ defined by the pure braid group $P_n$ on $n$ 
strands [the fundamental group of the space of ordered configurations of $n$ 
points in the plane] is generated by elements $x_{ik}, 1 \leq i < k \leq n$ 
subject to the relations
\[
[x_{ik},x_{st}] = 0 
\]
if $i,k,s,t$ are all distinct, and
\[
[x_{ik},x_{is}] = - [x_{ik},x_{ks}] 
\]
if $i,k,s$ are all distinct [32]. \bigskip

\noindent
The fundamental groups of the little disks operad define a (unital) operad 
$\{ P_n \}$ (the symmetric group action requires some care [47 \S 3]); in 
particular, cabling 
\[
c_I : P_{r(I)} \times \prod_{i \in I} P_i \to P_{|I|}
\]
of pure braids defines composition operations which extend to homomorphisms
\[
\cee_I : \p_{r(I)} \times \prod_{i \in I} \p_i \to \p_{|I|}
\]
of Lie algebras, defining an operad $\{ \p_n \}$ in that category as well. 
The natural product of Lie algebras is the direct sum of underlying
vector spaces, so $\cee_I$ is a sum of two terms, the second defined by the 
juxtaposition operation
\[
\p_{i_1} \oplus \cdots \oplus \p_{i_r} \to \p_{i_1 + \dots + i_r} \;.
\]
The remaining information is contained in a less familiar homomorphism
\[
\cee^0_I : \p_{r(I)} \to \p_{|I|}
\]
defined on the first component by any partition $I$ of $|I|$ in $r$ parts. It 
is not hard to see that
\[
\cee^0_I(x_{st}) \equiv \sum_{p \in i_s, q \in i_t} x_{pq} + \dots \;;
\] 
it would be very useful to know more about this expansion \dots
\bigskip

\noindent
{\bf 4} Groups act on themselves by conjugation, and thus in general to have 
lots of (inner) automorphisms; Lie algebras act similarly on themselves, by 
their adjoint representations. It would also be useful to understand something
about the relations between the automorphisms of an operad in groups or Lie 
algebras (as monoids in a category of objects $\{ \Oh_*\}$ with 
$\{ \Sigma_*\}$-action, as in \S 2 above), and systems of inner automorphisms 
of the objects $\Oh_n$: thus the adjoint action of a system of elements 
$\phi_n \in \p_n$ defines an operad endomorphism if $\cee^0_I(\phi_*) = 
\phi_{|I|}$ for all partitions $I$. From some perspective, the classification 
of such endomorphisms is really part of the theory of symmetric functions. 
\bigskip

\noindent
This may be relevant, because the action $\GT$ on the completed braid groups 
is relatively close to inner. Drinfel'd [15 \S 4] describes elements of $\GT$ 
as pairs $(\lambda,f)$, where $\lambda$ is a scalar (ie, an element of the 
field of definition for the kind of Lie algebras we're working with: in our 
case, $\Q$, for simplicity), and $f$ lies in the commutator subgroup of the 
Malcev completion of the free group on two elements. The pairs $(\lambda,f)$ 
are subject to certain restrictions which I'll omit here, on the grounds that 
the corresponding conditions on the Lie algebra are spelled out below. It is 
useful to regard elements $f = f(a,b)$ of the free group as noncommutative 
functions of two parameters $a,b$; if $\sigma_i$ is a standard generator of 
the braid group, and
\[
y_i = \sigma_{i-1} \cdots \sigma_2 \sigma_1^2 \sigma_2 \dots \sigma_{i-1} 
\in P_n 
\]
for $2 \leq i \leq n$ (with $y_1 = 1$), then the action of $\GT$ on the braid 
group is defined by
\[
(\lambda,f)(\sigma_i) = f(\sigma_i^2,y_i) \sigma_i^\lambda 
f(\sigma_i^2,y_i)^{-1} \;;
\]
the omitted conditions imply that when $i=1$, this reduces to an exponential 
automorphism $\sigma_1 \mapsto \sigma_1^{\lambda}$ defined  in the Malcev 
completion. \bigskip 

\noindent
{\bf 5} Here are some technical details about the Lie algebra of $\GT$, 
reproduced from [15 \S 5]. Drinfel'd observes [remark before Prop. 5.5] that 
the scalar term $\lambda$ can be used to define a filtration on this Lie 
algebra, and he describes the associated graded object $\grt$. The following 
formalism is useful: $\fr$ is the free formal $\Q$-Lie algebra defined by 
power series in two noncommuting generators $A,B$: it is naturally filtered 
by total polynomial degree, with the free graded Lie algebra on two 
generators as associated graded object. \bigskip

\noindent
The algebra $\grt$ consists of series $\psi = \psi(A,B) \in \fr$ which are 
antisymmetric [$\psi(A,B) = - \psi(B,A)$] and in addition satisfy the relations
\[
\psi(C,A)  + \psi(B,C) + \psi(A,B) = 0
\]
and 
\[
[B,\psi(A,B)] + [C,\psi(A,C)]  = 0
\]
when $A + B + C = 0$, as well as a third relation asserting that
\[
\psi(x_{12},x_{23} + x_{24}) + \psi(x_{13} + x_{23},x_{34}) - \psi(x_{12} + 
x_{13},x_{24} + x_{34})
\]
equals
\[ 
\psi(x_{23},x_{34}) + \psi(x_{12},x_{23}) \;,
\]
assuming that the $x_{ik}$ satisfy the relations defining $\p$. The bracket in 
$\grt$ is 
\[
\langle \psi_1,\psi_2 \rangle = [\psi_1,\psi_2] + \partial_{\psi_2}(\psi_1) -
\partial_{\psi_1}(\psi_2) \;,
\]
where $\partial_{\psi}$ is the derivation of $\fr$ given by $\partial_{\psi}
(A) = [\psi,A], \; \partial_{\psi}(B) = 0$. \bigskip

\noindent
Drinfel'd shows [15 \S 5.6] that this Lie algebra is in fact isomorphic to 
the Lie algebra of $\GT$ [omitting the subalgebra corresponding to the 
scalars, used in this description to define the grading], but the isomorphism 
is defined inductively, so describing its action on the braid groups is not 
immediate. Nevertheless, Ihara (cf. [15 \S 6.3]) has shown that for each odd 
$n > 1$ there are elements $\psi_n \in \grt$ such that
\[
\psi_n(A,B) \equiv \sum_{1 \leq m \leq n-1} \binom{n}{m} (\ad A)^{m-1} 
(\ad B)^{n-m-1} [A,B]
\]
modulo $[\fr',\fr']$ (where $\fr'$ is the derived Lie algebra of $\fr$].  
\bigskip

\noindent
It is conjectured that $\grt$ is free on these generators. Aside from [41], 
the cabling described in \S 3 doesn't seem to have been considered very 
closely, in this context. \bigskip

\newpage
                                                       
\bibliographystyle{amsplain}

\end{document}